\documentclass[11pt]{article}
\usepackage{amsfonts}
\usepackage{amscd}
\usepackage{amsmath}
\newlength{\dinwidth}
\newlength{\dinmargin}
\newtheorem{theorem}{Theorem}

\newtheorem{corollary}{Corollary}

\def\be{\begin{equation}}
\def\ee{\end{equation}}
\def\ben{\begin{displaymath}}
\def\een{\end{displaymath}}
\def\baa{\begin{eqnarray}}
\def\eaa{\end{eqnarray}}

\def\ba{\begin{array}}
\def\ea{\end{array}}

\makeatletter
\@addtoreset{equation}{section}
\makeatother

\def\f{\frac}
\def\la{\label}

\def\L{{\cal L}}

\def\p{\partial}
\def\B{{\bf B}}

\def\zbar{{\bar{z}}}

\begin{document}

\title{
Normalized Ricci flow on Riemann surfaces and determinants of
Laplacian }
\author{A. Kokotov\footnote{e-mail: alexey@mathstat.concordia.ca}   and
D. Korotkin
\footnote{e-mail: korotkin@mathstat.concordia.ca}}
\maketitle

\begin{center}
Department of Mathematics and Statistics, Concordia University\\
7141 Sherbrook West, Montreal H4B 1R6, Quebec,  Canada
\end{center}

\vskip2.0cm
{\bf Abstract.} In this note we give a simple proof of the fact that the determinant of Laplace operator
in smooth metric over compact Riemann surfaces of arbitrary genus $g$
 monotonously grows under the normalized Ricci flow. Together with
results of Hamilton that under the action of the normalized Ricci flow the 
smooth
metric tends asymptotically to metric of constant curvature for $g\geq 1$, this
leads to a simple proof of Osgood-Phillips-Sarnak
theorem stating that  within the class of smooth metrics with fixed
conformal class and fixed
 volume the determinant of Laplace operator
is maximal on metric of constant curvatute.

\vskip1.0cm

The normalized Ricci flow on compact  Riemann surface of arbitrary genus $g$
was  introduced by Hamilton \cite{ham}. Under the action of the
normalized Ricci flow the smooth metric $g_{ij}$ evolves according to
the following differential equation:
\be
\f{\p}{\p t} g_{ij} = (R_0-R) g_{ij}
\la{Ricci}
\ee
where $R$ is the scalar curvature, and $R_0$ is its average value. 
It was proved in \cite{ham} that solution of equation (\ref{Ricci})
exists for all times and that asymptotically, as $t\to\infty$, the metric converges to the metric of constant 
curvature for  $R_0\leq 0$, i.e. both for  $g=1$ (when  $R_0=0$)  and
$g\geq 2$ (when $R_0<0$). 
For $g=0$ the sufficient condition of convergency of the metric to
metric of constant curvature is the global condition $R > 0$ for
$t=0$.

The volume of the Riemann surface remains constant under the action of normalized Ricci flow (\ref{Ricci}).

In this note we give a simple  proof to  the following theorem
(although this statement in principle follows from results of
\cite{Sarnak}, the proof given here seems essentailly simpler):

\begin{theorem}
The determinant of Laplace operator in smooth metric over compact Riemann surface $\L$  monotonously 
increases under the action of Ricci flow (\ref{Ricci}). 
\end{theorem}

{\it Proof.} Write the metric in the diagonal form
$g_{ij}=\rho(z,\zbar)\delta_{ij}$; then the normalized Ricci flow 
takes the form 
\be
\f{\p}{\p t}\log\rho = R_0-R
\la{ricci1}
\ee
where 
\be
R=-\f{1}{\rho}(\log\rho)_{z\zbar}\;,
\la{curv}
\ee
\be
R_0=\f{\int_{\L}Rd\mu}{\int_{\L}\mu}
\la{R0}
\ee
where $\mu=\rho dz \wedge d\zbar$.
Consider  Alvarez-Polyakov formula \cite{Alvar,Polya} for variation of
$\det\Delta$ with respect to an arbitrary one-parametric infinitesimal
variation of  the metric in fixed conformal class:
 \be
\f{\p}{\p t}\left\{\f{\log\det\Delta}{\det\Im\B}\right\}
=\f{1}{12\pi}\int_{\L}(\log \rho)_t (-\rho R) dz\wedge d\zbar
\ee
which for the Ricci flow (\ref{ricci1}) 
equals
\be
\f{-1}{12}\int_{\L}(R_0-R) R d\mu = \f{-1}{12\int_{\L}d\mu}\left\{\left(\int_{\L} R d\mu\right)^2-
\int_{\L}d\mu \int_{\L} R^2 d\mu \right\}\geq 0
\la{cauchy}
\ee
according to Cauchy inequality $(\int_{\L} fg d\mu)^2 \leq (\int_{\L} f^2) (\int_{\L} g^2)$ for $f=1$, $g=R$.

Since the matrix of $b$-periods does not change under the Ricci flow, we get
\be
\f{\p}{\p t}\log\det\Delta\geq 0\;;
\ee
the equality takes place only if $R=const$ i.e. only for the metrics of constant curvature, 
which are stationary points of Ricci flow. 
\vskip0.5cm

The following corollary reproduces one of the results of the paper \cite{Sarnak}:

\begin{corollary}
On the space of smooth metrics on  Riemann  surfaces of  volume $1$
and genus $g\geq 1$
with fixed conformal structure the determinant 
has maximum on the metric of constant curvature ($0$ for $g=1$ and $-1$ for $d\geq 2$).
\end{corollary}

{\it Proof} follows from from Th.1 if we take into account results of
Hamilton \cite{ham} that in $g\geq 1$ asymptotically, as $t\to\infty$, any smooth
metric tends under the normalized Ricci flow  to the metric of
constant curvature (of the same volume).
\vskip0.5cm
For genus zero results of \cite{ham} imply a weaker version of Colollary 1, which states
that the determinat of Laplacian achieves its maximum on the metric of
constant curvature within the class of smooth metrics with fixed volume whose curvature
is everywhere positive.

We know a few  recent results  similar to Th.1. 
In \cite{Chen} it was proved that the
determinant of Laplace operator monotonously increases under Calabi flow \cite{Chen}.
In \cite{Ma} it was proved that the first eigenvalue of
Laplacian always 
increases under the action of (non-normalized) Ricci flow; however,
since the non-normalized Ricci flow does not preserve the volume, the
does not lead to a statement similar to Corollary 1 for the first
eigenvalue.

{\rm Remark (added December 2008)}. Recently Werner M\"uller informed us that the 
theorem proved in this note was proved in Section 3 of the
 earlier paper (\cite{Muller}.

\end{document}